\documentclass[12pt]{article}
\usepackage{amsmath,amsthm,amsfonts,amssymb,amscd,euscript}

\newlength{\defbaselineskip}
\newtheorem{theoremalpha}{Theorem}
\newtheorem{corollaryalpha}[theoremalpha]{Corollary}
\newtheorem{conjecturealpha}[theoremalpha]{Conjecture}

\theoremstyle{plain}
\newtheorem{thm}{Theorem}

\newtheorem{lem}[thm]{Lemma}

\newtheorem{prop}[thm]{Proposition}
\newtheorem{conj}[thm]{Conjecture}

\theoremstyle{definition}
\newtheorem{defn}[thm]{Definition}

\newtheorem{rem}[thm]{Remark}

\begin{document}

\begin{center}
{\large On the symmetric subscheme of Hilbert scheme of points

\vspace{.5cm}
 Kyungyong Lee}

\vspace{.3cm} Department of Mathematics, University of Michigan, Ann
Arbor, Michigan 48109, USA

kyungl@umich.edu
\end{center}



\begin{abstract}
We consider the Hilbert scheme $\text{Hilb}^{d+1}(\mathbb{C}^d)$ of
$(d+1)$ points in affine $d$-space $\mathbb{C}^d$ $(d\geq 3)$, which
includes the square of any maximal ideal. We describe equations for
the most symmetric affine open subscheme of
$\text{Hilb}^{d+1}(\mathbb{C}^d),$ in terms of Schur modules. In
addition we prove that $\text{Hilb}^{n}(\mathbb{C}^d)$ is reducible
for $n>d\geq 12$.
\end{abstract}




\section{Introduction}
Throughout these notes we work over $\mathbb{C}$. The maximal ideals
in a polynomial ring are very basic objects, and their deformations
are easy to understand. However very little is known about the
family of the ideals that can be deformed to the square of a maximal
ideal. Its existence and connectedness \cite{RH:conn} are well
known. Here we study its dimension.

We consider the Hilbert scheme $\text{Hilb}^{d+1}(\mathbb{C}^d)$ of
$(d+1)$ points in affine $d$-space $\mathbb{C}^d$, $d\geq 3$ (for
general introduction to the Hilbert schemes of points, see
\cite[\S18.4]{MS:comb}). It parametrizes  the ideals $I$ of colength
$(d+1)$ in $\mathbb{C}[\mathbf{x}]=\mathbb{C}[x_1, ... ,x_d]$. As
with any moduli problem, it is natural to ask whether
$\text{Hilb}^{d+1}(\mathbb{C}^d)$ is irreducible. It is already
interesting because $\text{Hilb}^{d+1}(\mathbb{C}^d)$ is irreducible
for $d\leq 3$ but reducible for $d\geq 12$.\footnote{Iarrobino
\cite{I:red} showed that $\text{Hilb}^{n}(\mathbb{C}^d)$ is
reducible for $d> 5$ and $n> (1+d)(1+d/4)$.}

\begin{theoremalpha}\label{reducible_example}
$\emph{Hilb}^{n}(\mathbb{C}^d)$ is reducible for $n>d\geq 12$.
\end{theoremalpha}

 Our main purpose is to describe equations for
the most symmetric affine open subscheme of
$\text{Hilb}^{d+1}(\mathbb{C}^d)$. We let $U\subset
\text{Hilb}^{d+1}(\mathbb{C}^d)$ denote the affine open subscheme
consisting of all ideals $I\in \text{Hilb}^{d+1}(\mathbb{C}^d)$ such
that $\{1,x_1, ... ,x_d\}$ is a $\mathbb{C}$-basis of
 $\mathbb{C}[\mathbf{x}]/I$. We will call $U$ \emph{the symmetric affine subscheme}.
 We note that the square of any maximal
 ideal in $\mathbb{C}[\mathbf{x}]$ belongs to the symmetric affine subscheme.

In these notes we give an elementary description of the coordinate ring of the symmetric
affine subscheme $U$. For a $\mathbb{C}$-vector space $V$ and a
partition $\lambda$, the module $\mathbb{S}_{\lambda}V$ is defined
by the Schur-Weyl construction. By abuse of notation, the quotient
ring given by the ideal generated by $\mathbb{S}_{\lambda}V$ in the
ring $\text{Sym}^{\bullet}(\mathbb{S}_{\mu}V)$ for some partitions
$\lambda$ and $\mu$ will be denoted by
$\frac{\text{Sym}^{\bullet}(\mathbb{S}_{\mu}V)}{<\mathbb{S}_{\lambda}V>}$.

\begin{theoremalpha}\label{mainthm1}
Let $d\geq 3$. Let $U$ be the symmetric affine open subscheme of
$\emph{Hilb}^{d+1}(\mathbb{C}^d)$. Then $U$ is isomorphic to
$$\mathbb{C}^d \times\emph{Spec}
\frac{\emph{Sym}^{\bullet}(\mathbb{S}_{(3,1,1,\cdots,1,0)}V)}{<\mathbb{S}_{(4,3,2,\cdots,2,1)}V>},$$where
$V$ is a $d$-dimensional $\mathbb{C}$-vector space,
$(3,1,1,\cdots,1,0)$ is a partition of $(d+1)$ and
$(4,3,2,\cdots,2,1)$ is of $(2d+2)$.
\end{theoremalpha}

Let us explain the notation more precisely. By Lemma \ref{comput},
there is an injective homomorphism
 $$j:\mathbb{S}_{(4,3,2,\cdots,2,1)}V \hookrightarrow \text{Sym}^{2}(\mathbb{S}_{(3,1,1,\cdots,1,0)}V)$$
of Schur modules. Then $j$ induces natural maps
$$\aligned &\mathbb{S}_{(4,3,2,\cdots,2,1)}V \otimes \text{Sym}^{r-2}(\mathbb{S}_{(3,1,1,\cdots,1,0)}V)\\
&\hookrightarrow
\text{Sym}^{2}(\mathbb{S}_{(3,1,1,\cdots,1,0)}V)\otimes
\text{Sym}^{r-2}(\mathbb{S}_{(3,1,1,\cdots,1,0)}V)\\
&\rightarrow \text{Sym}^{r}(\mathbb{S}_{(3,1,1,\cdots,1,0)}V),\text{
}\text{ }\text{ }\text{ }\text{ }\text{ }\text{ }\text{ }\text{
}\text{ }\text{ }\text{ }\text{ }\text{ }\text{ }\text{ }\text{
}\text{ }\text{ }\text{ }\text{ }\text{ }\text{ }\text{ }\text{
}\text{ }\text{ }\text{ }\text{ }\text{ }\text{ }\text{ }\text{
}\text{ }r \geq 2,\endaligned$$which define the quotient ring
$\frac{\text{Sym}^{\bullet}(\mathbb{S}_{(3,1,1,\cdots,1,0)}V)}{\mathbb{S}_{(4,3,2,\cdots,2,1)}V}$.

\begin{corollaryalpha}\label{maincor}
 Let $H(r)$ be the Hilbert function of
$\frac{\emph{Sym}^{\bullet}(\mathbb{S}_{(3,1,1,\cdots,1,0)}V)}{\mathbb{S}_{(4,3,2,\cdots,2,1)}V}$.
Let $$\emph{Sym}^{r}(\mathbb{S}_{(3,1,1,\cdots,1,0)}V)=\bigoplus_{|\lambda|=r(d+1)} \mathbb{S}_\lambda^{\text{ }\oplus m_\lambda} ,$$ where $m_\lambda\in \mathbb{Z}_{\geq 0}$ and $\lambda=(\lambda_1, \lambda_2, \cdots, \lambda_d)$ is a partition of $r(d+1)$, i.e., $\sum_{i=1}^d \lambda_i =r(d+1)$ and $\lambda_1\geq \lambda_2\geq \cdots\geq \lambda_d$.  Then
\begin{equation}\label{hilbfnbound}
H(r)\geq\sum_{\substack{|\lambda|=r(d+1) \\ \lambda_{d-k} +\cdots + \lambda_d\leq rk}} m_\lambda (\emph{dim}_{\mathbb{C}}  \mathbb{S}_\lambda),
\end{equation}
for any $r\geq 2$ and any $k=0,...,d-1$.
\end{corollaryalpha}

Corollary ~\ref{maincor} is an elementary consequence of the combinatorial Littlewood-Richardson rule(for example, see \cite[Appendix]{FH:repre}). In fact any  $\mathbb{S}_\lambda$ appearing in the decomposition of  $\mathbb{S}_{(4,3,2,\cdots,2,1)}V \otimes (\mathbb{S}_{(3,1,1,\cdots,1,0)}V)^{\otimes (r-2)}$ satisfies $\lambda_{d-k} +\cdots + \lambda_d\geq rk+1$, for any $r\geq 2$ and any $k=0,...,d-1$.

It is tedious but entirely possible to compute the right hand side
of (\ref{hilbfnbound}) for small $r$. These computations suggest
that the Hilbert function $H(r)$ grows faster than
$\mathcal{O}\left(r^{k {{d-k} \choose 2}}\right)$ for any
$k=0,...,d-1$. So Corollary ~\ref{maincor} suggests that, for
sufficiently large $d$, the symmetric open subscheme $U$ of
$\text{Hilb}^{d+1}(\mathbb{C}^d)$ has dimension greater than
$d(d+1)$, which implies that $\text{Hilb}^{d+1}(\mathbb{C}^d)$ is
reducible. To prove Theorem ~\ref{reducible_example}, we actually
find large dimensional families of ideals in a very explicit way.

\noindent\emph{Acknowledgements.} I am very much indebted to Boris
Shekhtman from whom I learned a great deal about ideal projectors.
Especially Lemma ~\ref{idealP} is due to him. I would like to thank
Carl de Boor for bringing my attention to ideal projectors, and
Steve Kleiman, Rob Lazarsfeld, Ezra Miller, John Stembridge, Dave
Anderson for valuable discussions.


\section{Proof of Theorem ~\ref{mainthm1}}

To ease notations and references, we introduce the notion of ideal
projectors(cf. \cite{B:algebra}, \cite{dB:ideal}, \cite{dB:limits},
\cite{S:limit}).

\begin{defn}
(cf. \cite{B:algebra}) A linear idempotent map $P$ on
$\mathbb{C}[\mathbf{x}]$ is called an \textbf{ideal projector} if
$\text{ker} P$ is an ideal in $\mathbb{C}[\mathbf{x}]$.
\end{defn}

We will use \emph{de Boor's formula}:

\begin{thm}\label{dB} \emph{(\cite{dB:ideal}, de Boor)} A linear mapping $P : \mathbb{C}[\mathbf{x}] \rightarrow \mathbb{C}[\mathbf{x}] $ is an
ideal projector if and only if the equality
\begin{equation}\label{deBoor}
P(gh) = P(gP(h))
\end{equation}
 holds for all $g, h \in \mathbb{C}[\mathbf{x}]$.
\end{thm}

Let $\mathcal{P}$ be the space of ideal projectors onto span
$\{1,x_1,...,x_d\}$, in other words,
$$ \mathcal{P}:=\{P : \text{ideal projector } | \text{ ker} P \in
U\}.
$$ The space $\mathcal{P}$ is isomorphic to the symmetric affine
subscheme $U$\cite[p3]{S:bivideal}. For the sake of simplicity, we
prefer to work on $\mathcal{P}$ in place of $U$.

First we consider the natural embedding of $\mathcal{P}$. Gustavsen,
Laksov and Skjelnes \cite{GLS:elem} gave more general description of
open affine coverings of Hilbert schemes of points.

\begin{lem}
The space $\mathcal{P}$ can be embedded into
$\mathbb{C}^{(d+1){{d+1}\choose 2}}$.
\end{lem}
\begin{proof}
For each ideal projector $P\in \mathcal{P}$ and each pair $(i,j)$,
$1\leq i,j\leq d$, there is a collection $p_{0,ij},p_{1,ij}, \cdots
,p_{d,ij}$ of complex numbers such that
$$P(x_i x_j) = p_{0,ij}+\sum_{m=1}^d p_{m,ij}x_m.$$
As $(i,j)$ varies over $1\leq i,j\leq d$, each ideal projector $P\in
\mathcal{P}$ gives rise to a collection $p_{0,ij}, p_{k, st}$
$(1\leq i,j,k,s,t \leq d)$ of complex numbers. Of course
$p_{0,ij}=p_{0,ji}$ and $p_{k,st}=p_{k, ts}$. So we have a map $f:
\mathcal{P} \rightarrow\mathbb{C}^{(d+1){{d+1}\choose
2}}=\frac{\mathbb{C}[p_{0,ij},\text{ }\text{ } p_{k, st}]_{1\leq
i,j,k,s,t \leq d}}{(p_{0,ij}-p_{0,ji},\text{ }\text{ }p_{k,st}-p_{k,
ts})}$.

Here we only show that $f$ is one-to-one. It is proved in
\cite{GLS:elem} that $f$ is in fact a scheme-theoretic embedding.

We will show that if $P_1, P_2\in\mathcal{P}$ and if
$f(P_1)=f(P_2)$, i.e. $P_1(x_i x_j)=P_2(x_i x_j)$ for every $(i,j)$,
$1\leq i,j\leq d$, then $P_1=P_2$. Since $P_1$ and $P_2$ are linear
maps, it is enough to check that
$P_1(x_{i_1}...x_{i_r})=P_2(x_{i_1}...x_{i_r})$ for any monomial
$x_{i_1}...x_{i_r}$. This follows from de Boor's formula
(\ref{deBoor}):
$$\aligned P_1(x_{i_1}...x_{i_r})&=P_1(x_{i_1}P_1(x_{i_2}\cdots
P_1(x_{i_{r-1}}x_{i_r})\cdots))\\
&=P_2(x_{i_1}P_2(x_{i_2}\cdots
P_2(x_{i_{r-1}}x_{i_r})\cdots))=P_2(x_{i_1}...x_{i_r}),
\endaligned$$where we have used the property that $P(g)$ is a linear
combination of $1,x_1,\dots,x_d$ for any $g \in
\mathbb{C}[\mathbf{x}]$.
\end{proof}

Next we describe the ideal defining $\mathcal{P}$ in
$$\frac{\mathbb{C}[p_{0,ij},\text{ }\text{ } p_{k, st}]_{1\leq
i,j,k,s,t \leq d}}{(p_{0,ij}-p_{0,ji},\text{ }\text{ }p_{k,st}-p_{k,
ts})}=:R,$$where we keep the notations in the above proof.  Let
$I_{\mathcal{P}}$ denote the ideal.

\begin{lem}\label{idealP}
Let $C(a;j,(i,k))\in R$ denote the coefficient of $x_a$ in
$$P(x_k P(x_i x_j)) - P(x_i P(x_k x_j)) \in R[x_1, \cdots,
x_d].$$ Then $I_{\mathcal{P}}$ is generated by $C(a;j,(i,k))$'s
$(0\leq a \leq d$,$\text{ }$ $1 \leq i,j,k\leq d)$. (We regard an
element in $R[x_1, \cdots, x_d]_0 \cong R$ as a coefficient of
$x_0$.)
\end{lem}

 For example, if $a\neq j,i,k$ then
$$C(a;j,(i,k))=\sum_{m=1}^d (p_{m,ij}p_{a,km}-p_{m,kj}p_{a,im}).$$
If $a=k$ then
$$C(k;j,(i,k))=p_{0,ij}+\sum_{m=1}^d (p_{m,ij}p_{k,km}-p_{m,kj}p_{k,im}).$$

\begin{proof}[Proof of Lemma ~\ref{idealP}]
The de Boor's formula (\ref{deBoor}) implies that $I_{\mathcal{P}}$
is generated by coefficients of $x_a$'s ($1\leq a \leq d$) in
$P(gP(h))-P(hP(g))$ (all $g,h\in \mathbb{C}[\mathbf{x}]$). But any
$P(gP(h))-P(hP(g))$ can be generated by $P(x_k P(x_i x_j)) - P(x_i
P(x_k x_j))$'s.
\end{proof}

We note that $C(a;j,(i,k))+C(a;j,(k,i))=0$ so from now on we
identify $C(a;j,(i,k))$ with $-C(a;j,(k,i))$.

\begin{lem}
In fact, $I_{\mathcal{P}}$ is generated by $C(a;j,(i,k))$'s $(1\leq
a \leq d$,$\text{ }$ $1 \leq i,j,k\leq d)$.
\end{lem}
\begin{proof}
It is enough to prove that for any $1 \leq i,j,k\leq d$, the
polynomial $C(0;j,(i,k))$ is generated by $C(a;b,(e,f))$'s $(1 \leq
a,b,e,f\leq d)$. Fix any $u$, $1 \leq u\leq d$. Then we have
$$\aligned C(0;j,&(i,k))
=\sum_{m=1}^d
(p_{m,ij}p_{0,km}-p_{m,kj}p_{0,im})\\
=-\sum_{m=1}^d \Big(&p_{m,ij} \sum_{t=1}^d
(p_{t,km}p_{u,tu}-p_{t,ku}p_{u,tm})
 -p_{m,kj}\sum_{t=1}^d(p_{t,im}p_{u,tu}-p_{t,iu}p_{u,tm})\Big)\\
+\sum_{m=1}^d \Big(&p_{m,ij} C(u;k,(m,u))
 -p_{m,kj}C(u;i,(m,u))\Big)
  \endaligned$$
 $$\aligned
 =-\sum_{t=1}^d \Big(&p_{u,tu}\sum_{m=1}^d (p_{m,ij}p_{t,km}-p_{m,kj}p_{t,im})\\
 -&p_{t,ku}\sum_{m=1}^d(p_{m,ij}p_{u,tm}-p_{m,it}p_{u,jm})
 +p_{t,iu}\sum_{m=1}^d(p_{m,kj}p_{u,tm}-p_{m,kt}p_{u,jm})\Big)\\
 +\sum_{m=1}^d
 &p_{u,jm}\sum_{t=1}^d(p_{t,ku}p_{m,it}-p_{t,iu}p_{m,kt})\\
 +\sum_{m=1}^d \Big(&p_{m,ij} C(u;k,(m,u))
 -p_{m,kj}C(u;i,(m,u))\Big)
 \endaligned$$
 $$\aligned
 =-\sum_{t=1}^d \Big(&p_{u,tu}C(t;j,(i,k))
 -p_{t,ku}C(u;i,(j,t))
 +p_{t,iu}C(u;k,(j,t))\Big)\text{ }\text{ }\text{ }\text{ }\text{ }\text{ }\text{ }\text{ }\text{ }\text{ }\text{ }\text{ }\text{ }\\
 +\sum_{m=1}^d
 &p_{u,jm}C(m;u,(k,i))\\
 +\sum_{m=1}^d \Big(&p_{m,ij} C(u;k,(m,u))
 -p_{m,kj}C(u;i,(m,u))\Big).
 \endaligned$$
\end{proof}

So the set of generators of $I_{\mathcal{P}}$ is
$$\{C(a;j,(i,k))\text{ } | \text{ }1 \leq a,i,j,k\leq d)\}.$$
We associate to this a representation of $GL(V)$.
\begin{prop}\label{vector_space}
The $\mathbb{C}$-vector space $W$ of generators
$$\frac{<C(a;j,(i,k))\text{ } | \text{ }1 \leq a,i,j,k\leq d)>}{C(a;j,(i,k))+C(a;j,(k,i))}$$
is canonically isomorphic to
$$
\mathbb{S}_{(3,2,1,\cdots,1,0)}V \bigoplus
\mathbb{S}_{(3,1,1,\cdots,1,1)}V$$ as $\mathbb{C}$-vector spaces,
where $V$ is a $d$-dimensional vector space and $
\mathbb{S}_{(3,2,1,\cdots,1,0)}$ $($resp.
$\mathbb{S}_{(3,1,1,\cdots,1,1)})$ is the Schur functor
corresponding to the partition $(3,2,1,\cdots,1,0)$ $($resp.
$(3,1,1,\cdots,1,1))$ of $(d+2)$.
\end{prop}
\begin{proof}
Let $V=\bigoplus_{i=1}^d \mathbb{C}v_i$. Define
$$\varphi:W \longrightarrow \bigwedge^{d-1}V \otimes V
\otimes \bigwedge^2 V$$ by
$$\varphi:C(a;j,(i,k)) \mapsto (-1)^a(v_1\wedge \cdots \wedge
\hat{v_a}\wedge \cdots \wedge v_d)\otimes v_j\otimes(v_i \wedge
v_k). $$ Then it is clear that $\varphi$ is injective.

By Littlewood-Richardson rule, we have
$$\aligned&\bigwedge^{d-1}V \otimes V \otimes \bigwedge^2 V\\
 &\cong \mathbb{S}_{(1,1,1,\cdots,1,0)}V \otimes V \otimes \mathbb{S}_{(1,1,0,\cdots,0,0)}
 V\\
 &\cong
\mathbb{S}_{(3,2,1,\cdots,1,0)}V \bigoplus
\mathbb{S}_{(3,1,1,\cdots,1,1)}V
 \bigoplus (\mathbb{S}_{(2,2,1,\cdots,1,1)}V)^{\oplus
2}\bigoplus \mathbb{S}_{(2,2,2,1,\cdots,1,0)}V\\
&\cong \mathbb{S}_{(3,2,1,\cdots,1,0)}V \bigoplus
\mathbb{S}_{(3,1,1,\cdots,1,1)}V
 \bigoplus\bigwedge^{d}V \otimes \bigwedge^2 V \bigoplus\bigwedge^{d-1}V \otimes
\bigwedge^3 V,\endaligned$$where each partition is of $(d+2)$. We
will show that the images of $W$ under $\varphi$ lie neither on
$\bigwedge^{d}V \otimes \bigwedge^2 V$ nor $\bigwedge^{d-1}V \otimes
\bigwedge^3 V$.

Since
$$
\sum_{j=1}^d (-1)^j(v_1\wedge \cdots \wedge \hat{v_j}\wedge \cdots
\wedge v_d)\otimes v_j\otimes(v_i \wedge v_k),\text{ }\text{ }\text{
}1\leq i <k\leq d,
$$
generate $\bigwedge^{d}V \otimes \bigwedge^2 V$, we need to show
that
\begin{equation}\label{eq1}
\sum_{j=1}^d C(j;j,(i,k)) =0.
\end{equation} But this is elementary
because
$$\sum_{j=1}^d C(j;j,(i,k))=\sum_{j=1}^d\sum_{m=1}^d (p_{m,ij}p_{j,km}-p_{m,kj}p_{j,im})=0.$$

Since
$$\aligned
&(v_1\wedge \cdots \wedge \hat{v_a}\wedge \cdots
\wedge v_d)\otimes v_j\otimes(v_i \wedge v_k)\\
&+(v_1\wedge \cdots \wedge \hat{v_a}\wedge \cdots
\wedge v_d)\otimes v_k\otimes(v_j \wedge v_i)\\
&+(v_1\wedge \cdots \wedge \hat{v_a}\wedge \cdots \wedge v_d)\otimes
v_i\otimes(v_k \wedge v_j),\text{ }\text{ }\text{ }1\leq a\leq
d,\text{ } 1\leq j<i<k\leq d,
\endaligned$$
generate $\bigwedge^{d-1}V \otimes \bigwedge^3 V$, we need to show
that
\begin{equation}\label{eq2}
C(a;j,(i,k))+C(a;k,(j,i))+C(a;i,(k,j)) =0.
\end{equation}
But this is again elementary because
$$
\aligned
&\sum_{m=1}^d (p_{m,ij}p_{a,km}-p_{m,kj}p_{a,im})\\
&+\sum_{m=1}^d (p_{m,jk}p_{a,im}-p_{m,ik}p_{a,jm})\\
&+\sum_{m=1}^d (p_{m,ki}p_{a,jm}-p_{m,ji}p_{a,km})=0.
\endaligned
$$

Therefore $\varphi(W) \subset \mathbb{S}_{(3,2,1,\cdots,1,0)}V
\bigoplus \mathbb{S}_{(3,1,1,\cdots,1,1)}V$, in other words,
$$\varphi:W \longrightarrow \mathbb{S}_{(3,2,1,\cdots,1,0)}V
\bigoplus \mathbb{S}_{(3,1,1,\cdots,1,1)}V$$ is injective.

The next lemma completes the proof.
\end{proof}

\begin{lem}
$\varphi:W \longrightarrow \mathbb{S}_{(3,2,1,\cdots,1,0)}V
\bigoplus \mathbb{S}_{(3,1,1,\cdots,1,1)}V$ is surjective.
\end{lem}
\begin{proof}
It is enough to show that there are no other nontrivial
$\mathbb{C}$-linear relations among $C(a;j,(i,k))$'s than
$\mathbb{C}$-linear combinations of (\ref{eq1}) and (\ref{eq2}).

Suppose
\begin{equation}\label{eq3}
C(a;j,(i,k))+\sum_{u,b,e,f}
c_{u,b,(e,f)}C(u;b,(e,f)) =0,\text{ }\text{ }\text{ }\text{ }\text{
}\text{ }\text{ }\text{ }c_{u,b,(e,f)}\in \mathbb{C}.
\end{equation}

If $a\neq i,j,k$ then $C(a;j,(i,k))$ contains a term
$p_{m,ij}p_{a,km}$ and a term $p_{m,kj}p_{a,im}$.  The term
$p_{m,ij}p_{a,km}$ appears only in $C(a;j,(i,k))$ and $C(a;i,(k,j))$
among all $C(u;b,(e,f))$, $1\leq u,b,e,f\leq d$. Similarly the term
$p_{m,kj}p_{a,im}$ appears only in $C(a;j,(i,k))$ and
$C(a;k,(j,i))$. So the left hand side of (\ref{eq3}) must be a
nontrivial linear combination of (\ref{eq2}) and other relations.

Similarly even if $a=i,j,$ or $k$, each term in $C(a;j,(i,k))$
appears only in the ones involved in (\ref{eq1}) or (\ref{eq2}). To
get cancelation among these, the left hand side of (\ref{eq3}) must
contain (\ref{eq1}) or (\ref{eq2}). Repeating the argument,
(\ref{eq3}) becomes a linear combination of (\ref{eq1}) and
(\ref{eq2}).
\end{proof}

\begin{lem}\label{comput}
We have $$\bigwedge^{d-1}V \otimes \emph{Sym}^2 V\cong
\mathbb{S}_{(2,1,1,\cdots,1,1)}V\oplus
\mathbb{S}_{(3,1,1,\cdots,1,0)}V,$$and
$$\aligned
\emph{Sym}^{2}(\mathbb{S}_{(3,1,1,\cdots,1,0)}V)\cong
&\mathbb{S}_{(6,2,2,\cdots,2,0)}V \oplus
\mathbb{S}_{(5,3,2,\cdots,1,1)}V \oplus
\mathbb{S}_{(5,2,2,\cdots,2,1)}V\\
& \oplus \mathbb{S}_{(4,4,2,\cdots,2,0)}V \oplus
\mathbb{S}_{(4,3,2,\cdots,2,1)}V \oplus
\mathbb{S}_{(4,2,2,\cdots,2,2)}V. \endaligned$$ (If $d=3$ then
$\mathbb{S}_{(5,3,2,\cdots,1,1)}V$ does not appear.)
\end{lem}
\begin{proof}
The first isomorphism follows from the Littlewood-Richardson rule.
The second isomorphism can be calculated by
\cite[pp.124--128]{R:comb}.
\end{proof}

\begin{lem}\label{firstanaly}
There is an injective homomorphism
 $$j:\mathbb{S}_{(4,3,2,\cdots,2,1)}V \hookrightarrow \emph{Sym}^{2}(\bigwedge^{d-1}V \otimes \emph{Sym}^2
V)$$ such that $\mathcal{P}$ (hence the symmetric affine subscheme
$U$) is isomorphic to
$$\emph{Spec}
\frac{\emph{Sym}^{\bullet}(\bigwedge^{d-1}V \otimes \emph{Sym}^2
V)}{\mathbb{S}_{(4,3,2,\cdots,2,1)}V},$$ where $(4,3,2,\cdots,2,1)$
is a partition of $(2d+2)$.
\end{lem}
\begin{proof}
Consider a diagram
$$\begin{CD} \frac{\mathbb{C}[p'_{0,ij}, \text{ }p'_{k, st}]_{1\leq i,j,k,s,t
\leq d}}{(p'_{0,ij}-p'_{0,ji}, \text{ }p'_{k,st}-p'_{k,ts})} @<f<<
\frac{\mathbb{C}[p_{0,ij}, \text{ }p_{k, st}]_{1\leq
i,j,k,s,t \leq d}}{(p_{0,ij}-p_{0,ji}, \text{ }p_{k,st}-p_{k,ts})}=:R\\
@VVgV     @.\\ T:=\frac{\mathbb{C}[p'_{k, st}]_{1\leq k,s,t \leq
d}}{(p'_{k,st}-p'_{k,ts})}  @.
\end{CD}$$
where $g$ is the natural projection and  $f^{-1}$ is defined by
$$\aligned
&p'_{0,ij} \mapsto C(i+1; j,(i,i+1)),\text{ }\text{ }\text{ }\text{
}\text{ }\text{ }\text{ }\text{ }\text{ }1\leq i\leq
j\leq d,\\
&\text{ }\text{ }\text{ }\text{ }\text{ }\text{ }\text{ }\text{
}\text{ }\text{ }\text{
}\text{ }\text{ }\text{(if } i=d \text{ then }i+1:=1) \\
 &p'_{k,st} \mapsto p_{k,st},\text{ }\text{ }\text{ }\text{ }\text{ }\text{ }\text{ }\text{ }\text{ }\text{ }\text{ }\text{ }\text{ }\text{ }\text{ }\text{ }\text{ }\text{ }\text{ }\text{ }\text{ }\text{ }\text{ }\text{ }\text{ }\text{ }\text{ }\text{
}1\leq s\leq t\leq d,\text{ }1\leq k\leq d.\endaligned$$ In fact $f$
is an isomorphism because $p_{0,ij}$ is a linear term in
$$C(i+1;j,(i,i+1))=p_{0,ij}+\sum_{m=1}^d
(p_{m,ij}p_{(i+1),(i+1)m}-p_{m,(i+1)j}p_{(i+1),im}).$$

Since $C(i+1;j,(i,i+1))\in I_{\mathcal{P}}$, we have an induced
isomorphism
\begin{equation}\label{equ1}\frac{R}{I_{\mathcal{P}}}\cong\frac{T}{I_{\mathcal{P}}T},\end{equation}
where $I_{\mathcal{P}}T$ is the expansion of $I_{\mathcal{P}}$ to
$T$. We note that in this construction $C(i+1; j,(i,i+1))$ can be
replaced by any $C(k;j,(i,k))$ or $C(k;i,(j,k))$ ($k\neq i,j$),
because the resulting $I_{\mathcal{P}}T$ does not depend on the
choice $C(k;j,(i,k))$ or $C(k;i,(j,k))$. In fact this construction
is natural in the sense that we eliminate all the linear terms
appearing in $C(a;j,(i,k))$ so that the ideal $I_{\mathcal{P}}T$ is
generated by quadratic equations.

Since $p'_{0,ij}$ are eliminated under passing $g$, the direct
summand $\mathbb{S}_{(3,1,1,\cdots,1,1)}V(\cong \text{Sym}^2 V)$ in
$W$ is eliminated. Then, by Proposition~\ref{vector_space}, the
vector space of generators of $I_{\mathcal{P}}T$ is canonically
isomorphic to $\mathbb{S}_{(3,2,1,\cdots,1,0)}V$ hence to
$$\bigwedge^{d}V \otimes\mathbb{S}_{(3,2,1,\cdots,1,0)}V \cong
\mathbb{S}_{(4,3,2,\cdots,2,1)}V \subset
\text{Sym}^2(\bigwedge^{d-1}V \otimes \text{Sym}^2 V),$$where the
last containment follows from Lemma~\ref{comput}.

The isomorphism of rings $$T =\frac{\mathbb{C}[p'_{k, st}]_{1\leq
k,s,t \leq d}}{(p'_{k,st}-p'_{k,ts})}\cong
\text{Sym}^{\bullet}(\bigwedge^{d-1}V \otimes \text{Sym}^2 V)$$
naturally induces the isomorphism of quotient rings
\begin{equation}\label{equ2}\frac{T}{I_{\mathcal{P}}T} \cong \frac{\text{Sym}^{\bullet}(\bigwedge^{d-1}V \otimes \text{Sym}^2
V)}{\mathbb{S}_{(4,3,2,\cdots,2,1)}V}.\end{equation}Combining this
with\begin{picture}(5,5)\put(4,0){(}\end{picture}~\ref{equ1}) gives
the desired result.
\end{proof}

\begin{thm}\label{mainthm}
There is an injective homomorphism
 $$j:\mathbb{S}_{(4,3,2,\cdots,2,1)}V \hookrightarrow \emph{Sym}^{2}(\mathbb{S}_{(3,1,1,\cdots,1,0)}V)$$
 such that $\mathcal{P}$ (hence the symmetric affine subscheme
$U$) is isomorphic to
$$\mathbb{C}^d \times\emph{Spec}
\frac{\emph{Sym}^{\bullet}(\mathbb{S}_{(3,1,1,\cdots,1,0)}V)}{\mathbb{S}_{(4,3,2,\cdots,2,1)}V},$$where
$(3,1,1,\cdots,1,0)$ is a partition of $(d+1)$ and
$(4,3,2,\cdots,2,1)$ is of $(2d+2)$.
\end{thm}
\begin{proof}[Sketch of Proof]
Define an isomorphism of rings $$T=\frac{\mathbb{C}[p'_{k,
st}]_{1\leq k,s,t \leq d}}{(p'_{k,st}-p'_{k,ts})}
\overset{\cong}\longrightarrow \frac{\mathbb{C}[q_{k, st}]_{1\leq
k,s,t \leq
d}}{(q_{k,st}-q_{k,ts})}=:Q$$by\\
\begin{picture}(200,30)
\put(100,12){$p'_{k,st} \mapsto $}
\put(135,12){\Big{\{}}\put(145,23){$q_{k,sk} + q_{s,ss},\text{
}\text{ }\text{ }\text{ }\text{ }\text{ }\text{ }\text{ }\text{ if }
k=t$}\put(145,2){$q_{k,st},\text{ }\text{ }\text{ }\text{ }\text{
}\text{ }\text{ }\text{ }\text{ }\text{ }\text{ }\text{ }\text{
}\text{ }\text{ }\text{ }\text{ }\text{ }\text{if } k\neq s,t.$}
\end{picture}

As a matter of fact this is a natural isomorphism, because the
square of any maximal ideal in $\mathbb{C}[\mathbf{x}]$ satisfies
$p'_{k,sk}-\frac{1}{2}p'_{s,ss}=0$ ($k\neq s$), i.e. $q_{k,sk}=0$.
It is straightforward to check that no element in minimal generators
of $I_{\mathcal{P}}Q$ contains terms involving $q_{s,ss}$, $1\leq
s\leq d$. For example, if $a,i,j,k$ are distinct, then
$$\aligned C(a;j,(i,k))&=\sum_{m=1}^d (p_{m,ij}p_{a,km}-p_{m,kj}p_{a,im})\\
&=\sum_{m\neq a,j,i,k} (p_{m,ij}p_{a,km}-p_{m,kj}p_{a,im})\\
&\text{ }\text{
}+(p_{j,ij}p_{a,kj}-p_{j,kj}p_{a,ij})+(p_{a,ij}p_{a,ka}-p_{a,kj}p_{a,ia})\\
&\text{ }\text{
}+(p_{i,ij}p_{a,ki}-p_{i,kj}p_{a,ii})+(p_{k,ij}p_{a,kk}-p_{k,kj}p_{a,ik})\endaligned$$
becomes
$$\aligned
&\sum_{m\neq a,j,i,k} (q_{m,ij}q_{a,km}-q_{m,kj}q_{a,im})\\
&\text{ }\text{
}+((q_{j,ij}+q_{i,ii})q_{a,kj}-(q_{j,kj}+q_{k,kk})q_{a,ij})+(q_{a,ij}(q_{a,ka}+q_{k,kk})-q_{a,kj}(q_{a,ia}+q_{i,ii}))\\
&\text{ }\text{
}+((q_{i,ij}+q_{j,jj})q_{a,ki}-q_{i,kj}q_{a,ii})+(q_{k,ij}q_{a,kk}-(q_{k,kj}+q_{j,jj})q_{a,ik})\\
&=\sum_{m\neq a,j,i,k} (q_{m,ij}q_{a,km}-q_{m,kj}q_{a,im})\\
&\text{ }\text{
}+(q_{j,ij}q_{a,kj}-q_{j,kj}q_{a,ij})+(q_{a,ij}q_{a,ka}-q_{a,kj}q_{a,ia})\\
&\text{ }\text{
}+(q_{i,ij}q_{a,ki}-q_{i,kj}q_{a,ii})+(q_{k,ij}q_{a,kk}-q_{k,kj}q_{a,ik}),
\endaligned$$
in which no term involves $q_{s,ss}$, $1\leq s\leq d$.

Therefore we get $$\frac{T}{I_{\mathcal{P}}T} \cong
\frac{Q}{I_{\mathcal{P}}Q} \cong \mathbb{C}[q_{s, ss}]_{1\leq s\leq
d} \otimes_{\mathbb{C}} {\frac{\mathbb{C}[q_{k, st}]_{1\leq k,s,t
\leq d,\text{ }\text{ } k\neq s\text{ or }t\neq
s}}{(q_{k,st}-q_{k,ts})}}\Big{/}{I_{\mathcal{P}}Q}.$$

On the other hand, Lemma~\ref{comput} implies
$$\text{Sym}^{\bullet}(\bigwedge^{d-1}V \otimes \text{Sym}^2 V)\cong
\text{Sym}^{\bullet}(\mathbb{S}_{(2,1,1,\cdots,1,1)}V\oplus
\mathbb{S}_{(3,1,1,\cdots,1,0)}V).$$ We may identify the basis of
$\mathbb{S}_{(2,1,1,\cdots,1,1)}V$ with $\{q_{s,ss} | 1\leq s\leq
d\}$. So,
by\begin{picture}(5,5)\put(4,0){(}\end{picture}~\ref{equ2}), we have
$$\aligned\frac{T}{I_{\mathcal{P}}T} &\cong
\mathbb{C}[q_{s, ss}]_{1\leq s\leq d} \otimes_{\mathbb{C}}
{\frac{\mathbb{C}[q_{k, st}]_{1\leq k,s,t \leq d,\text{ }\text{ }
k\neq s\text{ or }t\neq
s}}{(q_{k,st}-q_{k,ts})}}\Big{/}{I_{\mathcal{P}}Q}\\&\cong
\text{Sym}^{\bullet}(\mathbb{S}_{(2,1,1,\cdots,1,1)}V)\otimes
\frac{\text{Sym}^{\bullet}(\mathbb{S}_{(3,1,1,\cdots,1,0)}V)}{\mathbb{S}_{(4,3,2,\cdots,2,1)}V}.\endaligned$$Combining
this with\begin{picture}(5,5)\put(4,0){(}\end{picture}~\ref{equ1})
gives the desired result.
\end{proof}


\section{Proof of Theorem \ref{reducible_example}}

The following lemma is elementary and well-known (for example, see
\cite[Theorem 18.32]{MS:comb}). For the convenience of the reader,
we include its proof here.
\begin{lem}\label{elemlem}
If $\emph{Hilb}^{d+1}(\mathbb{C}^d)$ is reducible then so is
$\emph{Hilb}^{n}(\mathbb{C}^d)$ for $n\geq d+1$.
\end{lem}
\begin{proof}
Let $H_n^d:=\text{Hilb}^{n}(\mathbb{C}^d)$ and $R_n^d \subset H_n^d$
denote the closure of the open set parametrizing radical ideals. It
is enough to show that $R_n^d \neq H_n^d$ implies $R_{n+1}^d \neq
H_{n+1}^d$. If $I \in H_n^d \setminus R_n^d$ and $P = (p_1, ...
,p_d) \in \mathbb{C}^d$ is not a zero of $I$, then $I \cap \langle
x_1-p_1,...,x_d-p_d \rangle$ is a point in $H_{n+1}^d \setminus
R_{n+1}^d$.
\end{proof}

\begin{proof}[Proof of Theorem \ref{reducible_example}]
Due to Lemma \ref{elemlem}, it suffices to show that
$\text{Hilb}^{d+1}(\mathbb{C}^d)$ is reducible for $d\geq 12$.
Actually we will prove that the most symmetric open affine subscheme
$U$ of $\text{Hilb}^{d+1}(\mathbb{C}^d)$ is reducible for $d\geq
12$.

The symmetric open affine subscheme $U$ of
$\text{Hilb}^{14}(\mathbb{C}^{13})$ is reducible. In fact it can be
obtained by modifying Iarrobino's and Shekhtman's constructions (\cite{I:red}, \cite{S:ideal}). To each
$36\times 5$ matrix $B$ over $\mathbb{C}$, we associate an ideal
$$I_B:=\left(
  \begin{array}{c}
    x_6^2 + b_{1,1}x_1+b_{1,2}x_2+\cdots+b_{1,5}x_5, \\
    x_6 x_7 + b_{2,1}x_1+b_{2,2}x_2+\cdots+b_{2,5}x_5, \\
    \vdots \\
    x_{13}^2 + b_{36,1}x_1+b_{36,2}x_2+\cdots+b_{36,5}x_5 \\
  \end{array}
\right)+(x_1,\cdots,x_5)\cap\mathfrak{m}^2 +\mathfrak{m}^3,
$$where $\mathfrak{m}$ denotes the maximal ideal corresponding to the
origin. Note that $\{1,x_1,...,x_{13}\}$ is a $\mathbb{C}$-basis of
$\mathbb{C}[x_1, ..., x_{13}]/I_B$. Then the dimension of the family
of $I_B$'s is $36\cdot5=180$. Applying affine translations, we get a
$193(=180+13)$-dimensional family of ideals. Since $193> 13\cdot
14$(=the dimension of the closure of the set of radical ideals), $U$ is
reducible.

This construction can be easily generalized to $d\geq 13$. So the
symmetric open affine subscheme of
$\text{Hilb}^{d+1}(\mathbb{C}^{d})$ is reducible for $d\geq 13$. For
the case of $d=12$, we use $36\times 4$ matrices and get a
$156(=36\cdot4+12)$-dimensional family $\mathcal{F}$ of ideals. Then
$\mathcal{F}$ has the same dimension ($156=12\cdot 13$) as the
closure of the set of radical ideals, but a general member in
$\mathcal{F}$ is not a radical ideal. So $\mathcal{F}$ is not
contained in the principal (radical) component.
\end{proof}

It would be interesting to find equations for the principal
(radical) component - the component containing radical ideals - of
its symmetric open affine subscheme. The ideal defining the
principal component contains the ideal generated by
$\mathbb{S}_{(4,3,2,\cdots,2,1)}V$.


\section{Questions and Examples}

It is well known (\cite{K:desing}) that if $d=3$ then $U$ is
isomorphic to a cone over the Pl\"{u}cker embedding of the
Grassmannian $G(2,6)$ with a three-dimensional vertex. So we have
the following :

\begin{rem}\label{secondcor}
Let $V$ be a $3$-dimensional vector space and $W$ a $6$-dimensional
vector space. Then
$$
\frac{\text{Sym}^{\bullet}(\mathbb{S}_{(3,1,0)}V)}{\mathbb{S}_{(4,3,1)}V}\cong
\frac{\text{Sym}^{\bullet}(\bigwedge^2 W)}{\bigwedge^4 W}.$$ \qed
\end{rem}

Generalizing Theorem ~\ref{mainthm1}, we have

\begin{conjecturealpha}\label{mainconj1}
Let $d\geq 3$ and $n={{d+m}\choose m}$ for some positive integer
$m$. Let $U_m\subset \emph{Hilb}^n(\mathbb{C}^d)$ denote the affine
open subscheme consisting of all ideals $I\in
\emph{Hilb}^n(\mathbb{C}^d)$ such that $\{\mathbf{x}^{\mathbf{u}}
\text{ }|\text{ }|\mathbf{u}|\leq m \}=\{1,x_1, ... ,x_d^m\}$ is a
$\mathbb{C}$-basis of
 $\mathbb{C}[\mathbf{x}]/I$. Then there are injective homomorphisms
 $$j_k:\mathbb{S}_{(3k+1,2k+1,2k,\cdots,2k,k)}V \hookrightarrow \emph{Sym}^{2}(\mathbb{S}_{(2k+1,k,k,\cdots,k,0)}V),
  \text{ }\text{ }\text{ }\text{ }\text{ }\text{ }1\leq k \leq m,$$
of Schur modules such that $U_m$ is isomorphic to the induced scheme
by $j_k$
$$\mathbb{C}^d \times \prod_{k=1}^m \emph{Spec}
\frac{\emph{Sym}^{\bullet}(\mathbb{S}_{(2k+1,k,k,\cdots,k,0)}V)}{\mathbb{S}_{(3k+1,2k+1,2k,\cdots,2k,k)}V},$$where
$V$ is a $d$-dimensional vector space over $\mathbb{C}$,
$(2k+1,k,k,\cdots,k,0)$ is a partition of $k(d+1)$ and
$(3k+1,2k+1,2k,\cdots,2k,k)$ is of $2k(d+1)$.
\end{conjecturealpha}

One of possible ways to obtain the dimension of $U$ might be to find
its Hilbert polynomial.

\begin{lem}
Let $H(r)$ be the Hilbert function of
$\frac{\emph{Sym}^{\bullet}(\mathbb{S}_{(3,1,1,\cdots,1,0)}V)}{\mathbb{S}_{(4,3,2,\cdots,2,1)}V}$.
Then
$$H(1)=d{{d+1} \choose 2}-d
$$ and
$$H(2)={{d{{d+1} \choose 2}-d+1}\choose 2}-\frac{d^2 (d^2 -4)}{3}.
$$
\end{lem}
\begin{proof}
We recall the fact that if $\mu=(\mu_1 \geq \cdots \geq \mu_d\geq
0)$ then $$\text{dim }\mathbb{S}_{\mu}=\prod_{1\leq i < j\leq
d}\frac{\mu_i-\mu_j+j-i}{j-i}$$(for example, see \cite[Theorem
6.3]{FH:repre}). It is straightforward to check that
$$\text{dim }\mathbb{S}_{(3,1,1,\cdots,1,0)}V=d{{d+1} \choose 2}-d$$ and
$$\text{dim }\mathbb{S}_{(4,3,2,\cdots,2,1)}V=\frac{d^2 (d^2 -4)}{3}.$$
\end{proof}

\begin{conj}
Let $\tilde{H}(r)$ be the Hilbert polynomial of
$\frac{\emph{Sym}^{\bullet}(\mathbb{S}_{(3,1,1,\cdots,1,0)}V)}{\mathbb{S}_{(4,3,2,\cdots,2,1)}V}$.
Then $\tilde{H}(r)=H(r)$ for every $r\geq 0$.
\end{conj}

\begin{rem}
The case of $d=3$ is well-known(cf. \cite{K:desing}). Actual
equations $C(a;j,(i,k))$ (or \cite[p242]{K:desing}) are relatively
simple so we can use the computer algebra system Macaulay 2. It
shows that if $d=3$ then
$$\tilde{H}(r)=H(r)=14{{r+8}\choose 8}-21{{r+7}\choose 7}+9{{r+6}\choose 6}-{{r+5}\choose 5}.$$
\end{rem}

\begin{conj}
If $d\geq 4$ then $$\aligned H(3)&=\emph{dim
}\emph{Sym}^{3}(\mathbb{S}_{(3,1,1,\cdots,1,0)}V)-\emph{dim
}\mathbb{S}_{(4,3,2,\cdots,2,1)}V \otimes
\mathbb{S}_{(3,1,1,\cdots,1,0)}V\\
&\text{ }+ \emph{dim }\left(\mathbb{S}_{(6,4,3,\cdots,3,2)}V \oplus
\mathbb{S}_{(5,4,4,3,\cdots,3,2)}V \oplus
\mathbb{S}_{(5,4,3,\cdots,3,3)}V\right)\\
&={{d{{d+1} \choose 2}-d+2}\choose 3}-\frac{d^2 (d^2
-4)\left(d{{d+1} \choose
2}-d\right)}{3}+\frac{d(d^2-4)(3d^2+1)}{12}.\endaligned$$
\end{conj}

\begin{rem}
The conjecture holds true for $d=3$. When $d=4$ or $5$, it coincides
with the results obtained by using actual equations $C(a;j(i,k))$
and running Macaulay 2. The three Schur modules
$\mathbb{S}_{(6,4,3,\cdots,3,2)}V$,
$\mathbb{S}_{(5,4,4,3,\cdots,3,2)}V$, and
$\mathbb{S}_{(5,4,3,\cdots,3,3)}V$ appear with multiplicity $>1$ in
the decomposition of $\mathbb{S}_{(4,3,2,\cdots,2,1)}V \otimes
\mathbb{S}_{(3,1,1,\cdots,1,0)}V$.
\end{rem}


\end{document}